\pdfoutput=1
\documentclass[a4paper,11pt]{article}

\usepackage[T1]{fontenc}
\usepackage{lmodern}
\usepackage{amsmath,amssymb}
\usepackage{booktabs}
\usepackage{array}
\usepackage{multirow}
\usepackage[margin=25mm]{geometry}
\usepackage{listings}
\usepackage[hidelinks]{hyperref}

\newcommand{\code}[1]{\texttt{#1}}
\newcommand{\OzI}{Ozaki scheme~I}
\newcommand{\OzII}{Ozaki scheme~II}

\title{Ozaki Scheme II Is Fast on CPUs Too:\\
  \large Multiple-Precision Matrix Multiplication on
  Intel AMX-INT8 and Arm SVE2-i8mm}
\author{Tomonori Kouya\thanks{Otemon Gakuin University.}}
\date{August 2026}

\begin{document}
\maketitle

\begin{abstract}
We implement \OzII{} (residue number system + Chinese remainder theorem),
which reduces multiple-precision dense matrix multiplication to a sequence of
low-precision, high-throughput integer or floating-point GEMMs, on CPUs.
Two backends are built on top of a shared CRT reconstruction stage:
(a) exact \code{INT8}$\times$\code{INT8}$\to$\code{INT32} tile products on
Intel AMX, and (b) binary64 DGEMM, the CPU construction of the original
\OzII{} paper.  On a two-socket Xeon Gold 6526Y (Emerald Rapids, 32 cores)
we evaluate significand precisions of 53--2048 bits and matrix dimensions
$N=256$--$8192$.  The results are always within 1~ulp of a high-precision
MPFR reference (essentially correctly rounded), while running up to
167$\times$ faster than a naive MPFR matrix product, up to 588$\times$
faster than BNCmatmul's Strassen multiplication, and 9--78$\times$ faster
than \OzI{} (FP64 slicing + OpenBLAS DGEMM).  The break-even point between
the two backends is $N\approx2048$: below it the binary64 backend wins
thanks to its smaller number of moduli, above it the AMX-INT8 backend wins
as the GEMMs dominate.

We further port the implementation to aarch64 (NVIDIA GB10:
Cortex-X925~$\times$10 + Cortex-A725~$\times$10).  Since this machine lacks
SME/SME2, the INT8 kernel uses the SMMLA matrix-product instruction of the
SVE2 i8mm extension.  We obtain an exact INT8 GEMM sustaining 6.5~Top/s and,
still within 1~ulp across all conditions, speedups of 14--89$\times$ over
BNCmatmul's \OzI{} (OpenBLAS-linked routine) and 6--19$\times$ over a
fairness-adjusted OzI-best variant (Section~\ref{sec:arm}).
The paper also includes a tutorial introduction to \OzII{}
(Section~\ref{sec:intro-oz2}) and a quantitative explanation of why this
seemingly GPU-specific technique is fast on CPUs as well
(Section~\ref{sec:whyfast}).
\end{abstract}

\section{Introduction}

A naive implementation of multiple-precision matrix multiplication on top of
a variable-precision library such as MPFR~\cite{mpfr2007} performs one
multiple-precision multiply, add, and rounding per inner-product term, i.e.\
$O(MNK)$ expensive operations.  The classical acceleration route is
operation-count reduction \`a la Strassen, which is known to be effective in
multiple precision as well~\cite{kouya2014,kouya2016}.  The family of
methods by Ozaki et al.\ instead replaces the multiple-precision arithmetic
itself: the inputs are decomposed into low-precision quantities and a small
number of calls to highly optimized low-precision GEMM do all the work.

\begin{description}
\item[\OzI{}~\cite{ozaki2012}] splits the significands into several
  low-precision slices and computes $AB=\sum_{i,j}A_iB_j$ as a sum of
  low-precision GEMMs.  The number of slices grows linearly with the
  precision $p$, so the number of GEMMs is $O(p^2)$
  ($L(L{+}1)/2$ with triangular truncation).  It has been applied to GPU
  Tensor Cores~\cite{mukunoki2020}, to INT8 matrix
  engines~\cite{ozimmu2024,uchino2025}, and to multiple-precision linear
  algebra~\cite{kouya2023}.
\item[\OzII{}~\cite{ozaki2-2025}] reduces the inputs modulo a set of small
  pairwise-coprime moduli $m_1,\dots,m_L$, runs one low-precision GEMM per
  modulus, and reconstructs the result with the Chinese remainder theorem
  (CRT).  The number of moduli is linear in the precision, so the number of
  GEMMs is $O(p)$.  The original paper demonstrates INT8 Tensor Cores on
  GPUs and FP64 arithmetic on CPUs, reaching more than 80~TFlop/s of
  binary64-equivalent performance; a follow-up extends the scheme to FP8
  units via quantization~\cite{ozfp8-2026}.
\end{description}

So far, \OzII{} has mostly been evaluated in the context of ``emulating
binary64 with GPU matrix engines.''  The claim of this paper is its title:
\textbf{\OzII{} is fast on CPUs too---and not only for binary64 emulation,
but as an arbitrary-precision (tens to thousands of bits) matrix product}.
Concretely, we implement \OzII{} as an arbitrary-precision GEMM with MPFR
input and output, providing two low-precision GEMM backends on top of a
single CRT reconstruction stage:
(a) exact \code{INT8}$\times$\code{INT8}$\to$\code{INT32} tile products on
Intel AMX, and (b) binary64 DGEMM (OpenBLAS~\cite{openblas}), the CPU
construction of the original paper.  We then port the same code base to
aarch64 (NVIDIA GB10) and evaluate an INT8 kernel built on the SMMLA
instruction of the SVE2 i8mm extension~\cite{arm-arm}
(Section~\ref{sec:arm}).  CRT reconstruction, MPFR I/O, and shared-exponent
handling are common to all backends.

Section~\ref{sec:intro-oz2} first explains \OzII{} itself at a tutorial
level, and Section~\ref{sec:whyfast} explains quantitatively why it is fast
on CPUs.  The precise formulation of the algorithm is given in
Section~\ref{sec:algo}, and the implementation from Section~\ref{sec:impl}
onward.

\section{A Tutorial on \OzII{}: Exact Matrix Products via Residue Systems}
\label{sec:intro-oz2}

This section explains the ideas behind \OzII{} assuming little background.
We follow how the seemingly special problem of ``multiple-precision matrix
multiplication'' turns into ``many matrix multiplications over tiny
integers''---the workload modern hardware is best at.

\subsection{The problem: multiple-precision multiplication is built from bricks}

Numbers with more than 53 bits of precision (binary64) cannot be handled
directly by hardware; they are represented in software as arrays of 64-bit
integers called limbs.  A 1024-bit significand, for instance, consists of
16 limbs.

Multiplying two multiple-precision numbers proceeds like schoolbook long
multiplication: all pairwise limb products are formed and carries are
propagated.  One 1024-bit multiplication is therefore $16\times16=256$
integer multiplications plus carry handling, plus normalization and
rounding to make it a floating-point operation.  In a general-purpose
library such as MPFR, this comes with variable-length data-structure
management, function calls, and branches, so that \textbf{one
multiplication costs hundreds of nanoseconds}
(in our measurements, one \code{mpfr\_mul}+\code{mpfr\_add} pair at
$p=212$ bits costs about 80~ns per core).
A matrix product performs $N^3$ multiply--adds; for $N=1024$ that is about
$10^9$ of them, i.e.\ several seconds even on 20 cores.

Meanwhile, the low-precision matrix engines of the very same CPUs (AMX,
SVE2-i8mm) execute trillions of 8-bit integer multiply--adds per second.
\OzII{} is a way to \emph{translate} multiple-precision multiplication into
work for these fast engines.

\subsection{Step 1: make the inputs integers (shared exponents)}

Residue arithmetic (below) does not apply to floating-point values, so the
inputs are first converted to integers.  The key observation is that entry
$(i,j)$ of the product depends only on row $i$ of $A$ and column $j$ of
$B$.  We can therefore pick one scale (a shared exponent) per row of $A$
and per column of $B$,
\[
  \hat a_{ik} = \mathrm{round}(a_{ik}\,2^{\,t-e_i}), \qquad
  \hat b_{kj} = \mathrm{round}(b_{kj}\,2^{\,t-f_j}),
\]
turning them into $t$-bit integers, so that
\[
  \sum_k a_{ik}b_{kj} \;\approx\; \Big(\sum_k \hat a_{ik}\hat b_{kj}\Big)
  \cdot 2^{\,e_i+f_j-2t} .
\]
The scales factor out of the inner product.  Hence \textbf{if the integer
inner product $\hat c_{ij}=\sum_k\hat a_{ik}\hat b_{kj}$ can be computed
exactly}, all that remains is to multiply back a power of two and round
once.  Errors can only enter at the initial $\mathrm{round}$ and the final
rounding, so choosing $t$ slightly larger than the target precision
(default $t=p+16+\lceil\log_2K\rceil$) yields nearly correctly rounded
results; our measurements stay below 1~ulp in every condition
(Section~\ref{sec:acc}).

However, $\hat c_{ij}$ can be as large as $K\cdot2^{2t}$---for
$p=212,\ K=4096$ that is an integer of \textbf{nearly 500 bits}.  Computing
it directly would bring us back to square one.  This is where residue
arithmetic enters.

\subsection{Step 2: split one big-integer computation into many small ones
  (residue number systems)}

Integers enjoy the following convenient property.  Choose pairwise coprime
moduli $m_1,\dots,m_L$ and let $P=m_1m_2\cdots m_L$; then
\[
  \text{any integer } z \text{ with } |z| < P/2 \text{ is uniquely
  recoverable from } (z\bmod m_1,\ \dots,\ z\bmod m_L)
\]
(the Chinese remainder theorem, CRT).  Moreover, addition and
multiplication commute with taking residues:
\[
  (x+y)\bmod m = ((x\bmod m)+(y\bmod m))\bmod m,\qquad
  (xy)\bmod m = ((x\bmod m)(y\bmod m))\bmod m .
\]
In other words, \textbf{instead of computing one 500-bit inner product we
may compute $L$ inner products over numbers of roughly 8 bits}.  The $L$
computations are completely independent, and each becomes one
low-precision GEMM.

\paragraph{A small worked example}
Take the inner product $\hat c = 9\cdot8+(-7)\cdot6+5\cdot(-4)=10$ and the
moduli $m_1{=}5,\ m_2{=}7,\ m_3{=}9$ ($P=315$).  Mapping every element to
its representative of smallest absolute value (the \emph{centered residue},
$|r|\le m/2$):
\[
\begin{array}{c|ccc|ccc|c}
 m & 9 & -7 & 5 & 8 & 6 & -4 & \text{residue inner product} \\ \hline
 5 & -1 & -2 & 0 & -2 & 1 & 1 & (-1)(-2)+(-2)(1)+0\cdot1 = 0 \\
 7 & 2 & 0 & -2 & 1 & -1 & 3 & 2\cdot1+0\cdot(-1)+(-2)\cdot3 = -4\equiv 3 \\
 9 & 0 & 2 & -4 & -1 & -3 & -4 & 0\cdot(-1)+2\cdot(-3)+(-4)(-4) = 10\equiv 1 \\
\end{array}
\]
From the residue triple $(0,3,1)$, the unique integer with $|z|<157$
satisfying $z\equiv0\ (5),\ z\equiv3\ (7),\ z\equiv1\ (9)$ is $z=10$: the
correct inner product is recovered.  Note that \textbf{every intermediate
number had one or two digits}.  Had the inputs been 500-bit numbers instead
of $-7$ and $9$, the arithmetic after taking residues would have involved
numbers of exactly the same small size.

\subsection{Step 3: when is a low-precision GEMM \emph{exact}?}

When the residue matrices $A^{(\ell)}=\hat A\bmod m_\ell$ and
$B^{(\ell)}=\hat B\bmod m_\ell$ are multiplied on a low-precision engine,
the result is mathematically exact as long as \textbf{no rounding and no
overflow can occur}.  This is what bounds the moduli:

\begin{itemize}
\item \textbf{INT8 engines (AMX / SMMLA)}:
  for $m\le255$ the centered residues satisfy $|r|\le127$ and fit in a
  signed byte.  The multiply--adds accumulate in 32-bit integers, and
  integer addition has no rounding, so exactness holds as long as no
  overflow occurs: the worst case $K\cdot127^2<2^{31}$ gives
  $K\le133\,146$.
\item \textbf{binary64 DGEMM}: integers up to $2^{53}$ are represented
  exactly in binary64, and IEEE~754 additions and multiplications of such
  integers are exact while the results stay below $2^{53}$.  From the
  worst case $K\,(m/2)^2\le2^{53}$ we get $m\le\sqrt{2^{55}/K}$, about
  21.5 bits per modulus at $K=4096$.
\end{itemize}

So although we speak of ``low precision,'' \textbf{no approximation
whatsoever is made}: the low-precision engine is used as an exact
multiply--accumulate unit for small integers.  That is the heart of
\OzII{}; error exists only in the quantization of Step~1 and in the final
rounding.

\subsection{Step 4: CRT reconstruction}

Once the $L$ residue matrices
$C^{(\ell)}=A^{(\ell)}B^{(\ell)}\bmod m_\ell$ are available, each output
entry $\hat c_{ij}$ is reconstructed by the CRT and scaled back into an
MPFR value.  Reconstruction costs $O(L)$ per entry over $N^2$ entries and
is kept light by a formulation that avoids multiple-precision division
entirely (direct CRT, Section~\ref{sec:directcrt}).

\subsection{Contrast with \OzI{}: slicing}

The earlier \OzI{}~\cite{ozaki2012} uses \textbf{slices cut from the top of
the significand}, $\hat A=\sum_i A_i,\ \hat B=\sum_j B_j$, and computes
$\hat A\hat B=\sum_{i,j}A_iB_j$ as a sum of low-precision GEMMs.
This is natural, but \textbf{all pairs of slices} are needed, so with
$s\propto p$ slices the GEMM count grows as $O(p^2)$ (still $s(s{+}1)/2$
with triangular truncation).  In \OzII{} the moduli never interact
(each residue system is closed under arithmetic), so the GEMM count is
$O(p)$.  At $p=1024$ bits, \OzI{} needs 1431 DGEMMs whereas
\OzII{}/binary64 needs 99 (Section~\ref{sec:gemmcount}); the gap widens
with precision.  This is the main reason why the measured speed ratios in
this paper grow with precision.

\section{Why It Is Fast on CPUs}
\label{sec:whyfast}

\OzII{} was proposed mainly with INT8 Tensor Core GPUs in
mind~\cite{ozaki2-2025}, but nothing in the speedup mechanism is
GPU-specific.  It boils down to three points.

\subsection{Reason 1: per-operation cost differs by 3--4 orders of magnitude}

Table~\ref{tab:opcost} compares the effective cost of ``one multiply--add''
on the two machines used in this paper.  An MPFR multiply--add involves
limb-array long multiplication, normalization, rounding, and
variable-length data management, costing tens of nanoseconds.
Matrix instructions, on the other hand, perform many multiply--adds per
instruction: SMMLA does $2\times2\times8=32$ multiply--adds (64 integer
ops) per instruction, and AMX \code{TDPBSSD} does
$16\times16\times64=16384$.  The resulting throughput gap is four to five
orders of magnitude.

\begin{table}[htbp]
\centering
\caption{Measured effective throughput of multiply--add operations
  (all cores)}
\label{tab:opcost}
\small
\begin{tabular}{llrr}
\toprule
Operation & Machine & Effective rate & Per multiply--add \\
\midrule
MPFR multiply--add ($p=212$) & GB10, 20 threads & $2.5\times10^8$ /s & 4.0 ns \\
MPFR multiply--add ($p=212$) & Xeon, 32 threads & $\sim6\times10^8$ /s & 1.7 ns \\
OpenBLAS DGEMM & GB10 & 0.33 TFlop/s & 6.1 ps \\
OpenBLAS DGEMM & Xeon & 2.3 TFlop/s & 0.87 ps \\
INT8 GEMM (SMMLA) & GB10 & 6.5 Top/s & 0.31 ps \\
INT8 GEMM (AMX) & Xeon & 51.6 Top/s & 0.039 ps \\
\bottomrule
\end{tabular}
\end{table}

This comparison alone would be unfair: one MPFR multiply--add processes $p$
bits at once, whereas one INT8 multiply--add covers only 8 bits.  A fair
account must be made \emph{per bit of precision}---that is Reason~3.

\subsection{Reason 2: the GEMM shape unlocks arithmetic density}

Matrix multiplication performs $O(N^3)$ operations on $O(N^2)$ data, so
with blocking, the data can be kept in cache and reused many times.  This
is why optimized BLAS reaches peak performance, and lowering the precision
makes the structure \emph{even more} favorable: when an element shrinks
from 8 bytes (binary64) to 1 byte (int8), a cache of the same capacity
holds $8\times$ larger tiles, and the number of operations per byte of
memory traffic grows by $8\times$ as well.

A naive MPFR matrix product can exploit none of this.  Each element is a
variable-length limb array behind a pointer, so neither vectorization nor
tiling applies, and the code is scalar, branchy, and indirect.  In other
words, the naive method ``has the same $O(N^3)$ algorithm but uses well
under 0.1\% of the hardware's arithmetic capability.''
The preprocessing of \OzII{} (integerization, residue generation, data
reshuffling) is an investment that reshapes the data into the exact form
the arithmetic engines demand; its cost is $O(L\,N^2)$, one order below
the $O(L\,N^3)$ of the GEMMs.  As $N$ grows the investment is amortized
and the effective performance approaches that of the GEMM engine
(phase breakdowns: Tables~\ref{tab:phase} and~\ref{tab:armphase}).

\subsection{Reason 3: per-bit economics and the $O(p)$ GEMM count}

The total GEMM work of \OzII{} is $2LN^3$ operations with $L$ moduli.
One modulus carries $\log_2 m$ bits of information (INT8: about 8 bits;
binary64: about 21.5 bits), so for a required width of $2t+\log_2K$ bits we
get $L\propto p$.  The relative merit of the backends reduces to
\[
  \text{time} \;\propto\; \frac{L}{\text{GEMM throughput } R}
  \;\propto\; \frac{1}{b\cdot R}
  \qquad (b=\text{bits per modulus}),
\]
i.e.\ a comparison of the \textbf{price per bit} $b\cdot R$.  On the GB10,
\[
  \frac{(b R)_{\mathrm{INT8}}}{(b R)_{\mathrm{FP64}}}
  = \frac{7.97\times6.5\ \mathrm{Top/s}}{21.5\times0.33\ \mathrm{TFlop/s}}
  \approx 7.3 ,
\]
so even after paying for $2.7\times$ more moduli, INT8 keeps a
$>7\times$ advantage (the same coefficient is about 8 on the Xeon).

In summary, \OzII{} is fast on CPUs because
(i) all multiplications are reduced to \textbf{rounding-free integer
GEMMs}, a form that exploits dense matrix instructions (AMX, SMMLA) at
full tilt;
(ii) GEMMs have high arithmetic density while the pre/post-processing is
one order cheaper; and
(iii) the number of GEMMs is linear in the precision (quadratic for
\OzI{}).
All three hold whenever a low-precision, high-throughput matrix multiplier
exists.  GPU Tensor Cores are merely one example; CPU AMX and SVE2-i8mm
(and, in the future, SME/SME2) qualify just as well.

\section{Algorithm}
\label{sec:algo}

\subsection{Overall flow}

We want $C=AB$ with $A\in\mathbb{R}^{M\times K}$,
$B\in\mathbb{R}^{K\times N}$, whose entries are MPFR floating-point numbers
of $p$ bits.  The computation has five stages:

\begin{quote}\ttfamily
MPFR input\\
\quad$\downarrow$ (1) shared-exponent integer significands ($t$ bits)\\
integer matrices $\hat A,\hat B$\\
\quad$\downarrow$ (2) centered residues per modulus\\
residue matrices $A^{(\ell)},B^{(\ell)}$\\
\quad$\downarrow$ (3) exact low-precision GEMM ($\times L$)\\
$C^{(\ell)}=A^{(\ell)}B^{(\ell)} \bmod m_\ell$\\
\quad$\downarrow$ (4) direct CRT reconstruction\\
signed integers $\hat c_{ij}$\\
\quad$\downarrow$ (5) rescale and round\\
MPFR output
\end{quote}

\subsection{Integerization with shared exponents}

For row $i$ of $A$, let $e_i=\max_k\mathrm{exp}(a_{ik})$ over its nonzero
entries (MPFR's convention is $|x|<2^{\mathrm{exp}(x)}$).  With an internal
significand width of $t$ bits,
\begin{equation}
  \hat a_{ik} = \mathrm{round}\!\left(a_{ik}\,2^{\,t-e_i}\right),
  \qquad |\hat a_{ik}|\le 2^{t},
\end{equation}
and similarly for column $j$ of $B$ with
$f_j=\max_k\mathrm{exp}(b_{kj})$ and
$\hat b_{kj}=\mathrm{round}(b_{kj}2^{\,t-f_j})$.  Then
\begin{equation}
  \hat c_{ij}=\sum_{k}\hat a_{ik}\hat b_{kj}\in\mathbb{Z},
  \qquad
  c_{ij}= \hat c_{ij}\,2^{\,e_i+f_j-2t},
\end{equation}
so once $\hat c_{ij}$ is exact, only one final rounding remains.  The
quantization error is bounded by roughly $K\,2^{t}$, so by default
\begin{equation}
  t = p + g + \lceil\log_2 K\rceil
  \qquad (g:\ \text{guard bits, default 16}).
\end{equation}

\subsection{Unique CRT reconstruction and the choice of moduli}

Since $|\hat c_{ij}|\le K\,2^{2t}$, a unique signed (balanced) CRT
reconstruction requires
\begin{equation}
  P=\prod_{\ell=1}^{L} m_\ell \;>\; 2K\,2^{2t}
  \quad\Longleftrightarrow\quad
  \log_2 P > 2t+\log_2 K+1 .
  \label{eq:crtcond}
\end{equation}
Our implementation imposes the safer $\log_2P>2t+\log_2K+2$.

The upper bound on the moduli comes from the requirement that the
per-modulus GEMM run \emph{without rounding error}, which depends on the
backend.

\paragraph{(a) AMX-INT8 backend}
For a centered residue $|r|\le\lfloor m/2\rfloor$ to fit in one signed
byte, $m\le255$ suffices ($|r|\le127$).  Accumulation is in \code{INT32},
so exactness holds up to $K\cdot127^2<2^{31}$, i.e.\ $K\le133\,146$ (the
implementation caps $K$ at $133\,000$).  Crucially, this removes the
multi-limb decomposition required by FP4/FP8-type
implementations~\cite{ozfp8-2026}: \textbf{one modulus costs exactly one
AMX GEMM}.

The pairwise-coprime moduli with $m\le255$ (the 54 prime powers up to 255)
have a product of only about $2^{361.8}$, which by \eqref{eq:crtcond}
limits reconstruction to $t\lesssim174$ bits.  Beyond that, moduli up to
$m\le16191$ are added and their centered residues are split in balanced
base~128: $r=r_0+128\,r_1$ with $|r_0|\le64$, $|r_1|\le63$, computed with
\textbf{three GEMMs} per modulus via Karatsuba,
\begin{equation}
  \sum_k r^{(a)}_k r^{(b)}_k
  = P_0 + 128\,(P_s-P_0-P_2) + 16384\,P_2,
\end{equation}
\begin{equation}
  P_0=\textstyle\sum a_0b_0,\quad
  P_s=\textstyle\sum (a_0{+}a_1)(b_0{+}b_1),\quad
  P_2=\textstyle\sum a_1b_1,
\end{equation}
where $|a_0{+}a_1|\le127$ still fits in \code{INT8}.  One-limb moduli give
7.97 bits/GEMM at the top (6.70 on average over the 54), two-limb moduli
13.98 bits per 3 GEMMs $=4.66$ bits/GEMM, so the hybrid uses up all
one-limb moduli before adding two-limb ones.

\paragraph{(b) binary64 backend}
Following the CPU construction of the original paper~\cite{ozaki2-2025},
residues are represented exactly in binary64 and multiplied by an ordinary
DGEMM.  The inner product of centered residues obeys
$|{\sum}|\le q\,(m/2)^2$, so with $u=2^{-53}$,
\begin{equation}
  q\,m^2 \le 4u^{-1}=2^{55}
  \quad\Longleftrightarrow\quad
  m \le \sqrt{2^{55}/q}
  \label{eq:fp64bound}
\end{equation}
guarantees that the DGEMM result is an exact integer in binary64 ($q$ is
the inner-product length; here $q=K$).  For $K=4096$ this allows
$m\le2\,965\,821$, about 21.5 bits per modulus.

\paragraph{Required number of GEMMs}
Table~\ref{tab:moduli} lists the parameters of both backends at $K=4096$.
The binary64 backend needs fewer moduli ($2.7\times$ more bits per
modulus); the AMX backend has vastly faster individual GEMMs.  At $t=240$
($p=212$, quad-word) the binary64 backend uses 23 moduli / 23 DGEMMs,
closely matching the 22 moduli reported in the original paper.

\begin{table}[htbp]
\centering
\caption{Moduli configuration ($K=4096$, $t=p+16+\lceil\log_2K\rceil$)}
\label{tab:moduli}
\small
\begin{tabular}{rr rrr rrr}
\toprule
& & \multicolumn{3}{c}{AMX-INT8} & \multicolumn{3}{c}{binary64} \\
\cmidrule(lr){3-5}\cmidrule(lr){6-8}
$p$ & $t$ & moduli $L$ & GEMMs & $\log_2P$ & moduli $L$ & DGEMMs & $\log_2P$ \\
\midrule
   53 &   81 &  24 &  24 &  182 &   9 &   9 &  194 \\
  106 &  134 &  39 &  39 &  283 &  14 &  14 &  301 \\
  113 &  141 &  42 &  42 &  301 &  14 &  14 &  301 \\
  212 &  240 &  64 &  84 &  502 &  23 &  23 &  495 \\
  256 &  284 &  70 & 102 &  586 &  28 &  28 &  602 \\
  512 &  540 & 107 & 213 & 1102 &  51 &  51 & 1097 \\
 1024 & 1052 & 181 & 435 & 2131 &  99 &  99 & 2129 \\
 2048 & 2076 & 329 & 879 & 4174 & 194 & 194 & 4171 \\
\bottomrule
\end{tabular}
\end{table}

\subsection{Direct CRT without multiple-precision division}
\label{sec:directcrt}

Let $M_\ell = P/m_\ell$ and $u_\ell = \hat c_{ij}\bmod m_\ell$, and set
\begin{equation}
  v_\ell = u_\ell\,(M_\ell^{-1}\bmod m_\ell)\bmod m_\ell,
  \qquad
  S = \sum_{\ell=1}^{L} v_\ell M_\ell .
\end{equation}
Because $S/P=\sum_\ell v_\ell/m_\ell$ holds \emph{exactly},
\begin{equation}
  S = x_0 + \alpha P,
  \qquad
  \alpha=\left\lfloor \sum_{\ell} \frac{v_\ell}{m_\ell} \right\rfloor,
  \qquad 0\le x_0<P,
\end{equation}
and $\alpha$ is a small integer with $\alpha<L$.  We therefore evaluate
$\alpha$ in binary64, conservatively use $\alpha-1$, compute
$R=S-(\alpha-1)P\in[0,3P)$ with integer arithmetic, and bring it into
$[0,P)$ with at most two subtractions, obtaining the \textbf{exact $x_0$
without any multiple-precision division}.  Finally
\begin{equation}
  \hat c_{ij}=\begin{cases}
    x_0 & (x_0\le \lfloor P/2\rfloor)\\
    x_0-P & (\text{otherwise})
  \end{cases}
\end{equation}
gives the signed representative.  This reorganizes the ``FP64 quotient
estimation + integer correction'' of the original
paper~\cite{ozaki2-2025} into a form where the quotient is determined
exactly.

\section{Implementation}
\label{sec:impl}

\subsection{Structure}

About 2900 lines of C, depending on GMP 6.3.0~\cite{gmp} /
MPFR 4.2.2~\cite{mpfr2007} / OpenBLAS 0.3.33~\cite{openblas}.
There are two INT8 GEMM kernels---x86-64 (AMX) and aarch64
(SVE2-i8mm)---selected automatically by the \code{Makefile} via
\code{uname -m}.  They share pack-buffer dimensions and function
interfaces; only the internal layouts differ.

\begin{table}[htbp]
\centering
\caption{Source files}
\small
\begin{tabular}{lrl}
\toprule
File & Lines & Contents \\
\midrule
\code{src/oz2.h}          & 132 & public API \\
\code{src/oz2\_ctx.c}     & 241 & modulus selection, CRT tables (mpn and 32-bit digits), allocator \\
\code{src/oz2\_amx.c}     & 206 & exact AMX-INT8 GEMM ($2\times2$ tiles, L2 blocking; x86-64) \\
\code{src/oz2\_i8mm.c}    & 161 & exact SVE2-i8mm (SMMLA) GEMM ($8\times8$ blocks; aarch64) \\
\code{src/oz2\_split.c}   & 550 & MPFR $\to$ integer significands $\to$ centered residues $\to$ tile layout \\
\code{src/oz2\_crt.c}     & 195 & direct CRT (scalar mpn and AVX-512 batched) \\
\code{src/oz2\_gemm.c}    & 382 & driver (two backends, two parallel schedules) \\
\bottomrule
\end{tabular}
\end{table}

The API is as follows; the backend is selected via
\code{oz2\_opts.backend}.

\begin{lstlisting}[language=C]
oz2_opts o;
oz2_opts_default(&o);
o.backend    = OZ2_FP64;   /* or OZ2_AMX_INT8 (default) */
o.guard_bits = 16;         /* t = prec + guard + ceil(log2 K) */
oz2_gemm(M, N, K, A, lda, B, ldb, C, ldc, &o);
\end{lstlisting}

\subsection{The AMX-INT8 GEMM kernel}

AMX~\cite{intel-sdm} provides eight tile registers ($16\times64$ bytes
each); \code{TDPBSSD} performs $C\mathrel{+}=A\cdot B$ as
\code{INT8}$\times$\code{INT8}$\to$\code{INT32}.  The $B$ operand must be
in the VNNI layout, interleaving groups of 4 along $K$:
\begin{equation}
  B_{\mathrm{VNNI}}[\lfloor k/4\rfloor][4j+(k\bmod4)] = B[k][j].
\end{equation}

With only eight tiles, the maximal configuration is a $2\times2$
accumulator block ($32\times32$ \code{INT32} outputs) plus two $A$ tiles
and two $B$ tiles, giving a 1:1 ratio of tile loads to \code{TMUL}s.

The decisive factor is \textbf{locality of the tile loads}.
Table~\ref{tab:tileload} shows \code{TMUL} throughput as a function of
working-set size: loads from L3 or memory cost 6--8$\times$ more than from
L1 and collapse performance by a factor of six.

\begin{table}[htbp]
\centering
\caption{\code{TDPBSSD} throughput vs.\ working set (1 core)}
\label{tab:tileload}
\small
\begin{tabular}{lrrr}
\toprule
Working set & Size & ns / TMUL & Gop/s \\
\midrule
registers only (no loads) &     4 KiB & 5.82 & 5632 \\
L1-resident               &    32 KiB & 5.73 & 5722 \\
L2-resident               &  1024 KiB & 7.55 & 4337 \\
L3                        &  8192 KiB & 36.47 & 899 \\
memory                    & 65536 KiB & 44.94 & 729 \\
\bottomrule
\end{tabular}
\end{table}

The kernel therefore blocks $K$/$M$/$N$ so that one $A$ block and one $B$
block stay in L2 (defaults \code{KCB}=64 tiles, \code{MCB}=8,
\code{NCB}=32), and places the $N$-block loop \emph{outside} the OpenMP
work-sharing loop so that all threads stream the same $B$ block at the same
time, keeping it in the shared L3.  This alone raised the 32-thread
effective rate from 5.0 to 27~Top/s.

For small matrices with many moduli (e.g.\ $N=128$, 1024 bits, 435
planes), the per-GEMM OpenMP region and barriers dominate; the driver then
switches automatically to a schedule that \emph{parallelizes across
modulus planes} and runs each GEMM single-threaded.

\subsection{Residue generation and packing}

The significands extracted from MPFR are kept as base-$2^{16}$ digit
arrays ($A$ in order \code{[i][d][k]}, $B$ in \code{[k][d][j]}), and the
residue modulo $m$ is computed as
\begin{equation}
  \hat a \bmod m \;=\; \Big(\sum_{d} \mathrm{dig}_d \cdot (2^{16d}\bmod m)\Big) \bmod m ,
\end{equation}
with the digit loop outside and the element loop inside, so the element
loop auto-vectorizes (AVX-512); the reduction uses multiplication by a
double-precision reciprocal plus correction.  The VNNI conversion of $B$
processes four $k$-rows at a time with a $4\times16$-byte transpose
(six \code{punpck}-class instructions) and writes each 64-byte
destination---exactly one cache line---with a non-temporal store
(\code{vmovntdq}).

\subsection{AVX-512 batched CRT}

The digit direction of multiple-precision arithmetic cannot be vectorized
because of carries, but \emph{distinct output elements} are fully
independent.  We therefore adopt the layout ``32-bit digits, SoA, eight
output elements per vector'': $S$ is stored as 32-bit digits, one digit per
64-bit lane of eight lanes, processed as
\begin{equation}
  t_k = S_k + v\cdot M_{\ell,k} + \mathrm{carry}_k,\quad
  S_k' = t_k \bmod 2^{32},\quad
  \mathrm{carry}_{k+1}=\lfloor t_k/2^{32}\rfloor
\end{equation}
using \code{\_mm512\_mul\_epu32} and \code{\_mm512\_srli\_epi64}.
The evaluation of $\alpha$, the subtraction $S-\alpha P$, the correction
into $[0,P)$, and the signed conversion are all branch-free mask
operations.  A scalar GMP mpn version (\code{mpn\_addmul\_1},
\code{mpn\_submul\_1}) is kept as the reference path; both are verified
against exact integers over the whole range $t=40$--$1120$ bits for the
moduli sets of both backends.

\paragraph{Digit headroom}
$S=\sum_\ell v_\ell M_\ell$ can reach $P\cdot\sum_\ell m_\ell$, which
overflows the digit count of $P$ itself.  Two extra 32-bit digits (64
bits) of headroom are allocated.  This is the kind of bug that only
manifests when $\log_2P$ is close to $64\lceil\log_2P/64\rceil$; it was
caught by a dense parameter sweep in the unit tests.

\subsection{Effect of the main optimizations}

\begin{table}[htbp]
\centering
\caption{Effect of the main optimizations}
\small
\begin{tabular}{p{78mm}p{62mm}}
\toprule
Item & Effect \\
\midrule
L2 blocking of AMX tile loads ($B$ block shared by all threads)
 & 5.0 $\to$ 27 Top/s (32 threads, $N$=4096)\\
\code{MADV\_HUGEPAGE} (2 MiB pages)
 & 1.6--2.4$\times$ overall at $N$=4096; large-stride non-temporal stores
   devastate a 4 KiB-page DTLB \\
digit loop hoisted out of the residue kernel + AVX-512
 & pack 0.30 $\to$ 0.16 s ($N$=2048, 113 bits)\\
VNNI conversion of $B$ via $4\times16$ transpose + \code{vmovntdq}
 & another $\sim$30\,\% off pack \\
modulus loop moved inside the digit-array sweep
 & digit arrays read $1/L$ times \\
AVX-512 batched CRT
 & CRT 2.0--2.2$\times$ faster \\
\bottomrule
\end{tabular}
\end{table}

\section{Experimental Setup (x86-64)}

\begin{table}[htbp]
\centering
\caption{Evaluation environment (x86-64)}
\small
\begin{tabular}{ll}
\toprule
CPU & Intel Xeon Gold 6526Y (Emerald Rapids) $\times$ 2 sockets \\
    & 16 cores/socket, 32 cores total, max 3.9 GHz, AMX-INT8 / AVX-512 \\
Caches & L1d 48 KiB, L2 2 MiB/core, L3 37.5 MiB/socket \\
Memory & 1 TiB, 2 NUMA nodes \\
OS & Linux 6.8.0 (AMX permission via \code{ARCH\_REQ\_XCOMP\_PERM}) \\
Compiler & GCC 13.3.0, \code{-O3 -march=native -mamx-tile -mamx-int8 -fopenmp} \\
Libraries & GMP 6.3.0, MPFR 4.2.2, OpenBLAS 0.3.33 (SAPPHIRERAPIDS kernels) \\
Runtime env & \code{OMP\_NUM\_THREADS=32 OMP\_PROC\_BIND=close OMP\_PLACES=cores} \\
Baseline & BNCmatmul 0.24~\cite{bncmatmul,kouya2014} (AVX-512 + OpenMP build) \\
\bottomrule
\end{tabular}
\end{table}

Input matrices are uniform random in $[-1,1)$ (some tests widen the
exponent range by $\pm$spread).  Accuracy is reported as the maximum
relative error, in ulps of the target precision ($2^{-p}$), of 32 randomly
chosen output entries against references computed with $p+256$-bit inner
products.

\section{Results (x86-64)}

\subsection{Correctness}

For $p=53,106,113,128,200,256,512,1024$ bits and various $M,N,K$ (square
and non-square, including widened exponent ranges), the outputs
\textbf{agreed in every element} with a $p+160$-bit MPFR reference rounded
to $p$ bits (0~ulp).  A standalone CRT test sweeps $t=40$--$1120$ bits in
steps of 8 for the moduli sets of both backends and confirms that the
scalar and AVX-512 paths reconstruct randomly generated exact integers
identically.

\subsection{Standalone AMX-INT8 GEMM performance}

Exactness (full agreement with naive computation) was confirmed at all
sizes.  A single isolated GEMM reaches about 20~Top/s at $N=4096$; inside
the \OzII{} pipeline, executing many modulus planes back to back, the rate
is 39--43~Top/s at $N=K=4096$ and up to \textbf{51.6~Top/s} at
$N=K=8192$.  Relative to the per-core L1-resident peak (5.72~Top/s,
Table~\ref{tab:tileload}) this is 23--28\% efficiency on 32 cores; the
remainder is mostly tile loads from L2/L3 and $C$-tile stores.

\subsection{Phase breakdown}

Table~\ref{tab:phase} shows the phase breakdown at $N=4096$.  In the AMX
backend the GEMMs account for only 13--17\% of the total; the dominant
phase is the $O(L\,N^2)$ residue conversion (pack).  In the binary64
backend, conversely, the GEMMs take 61--66\%.

\begin{table}[htbp]
\centering
\caption{Phase breakdown ($N=M=K=4096$, seconds)}
\label{tab:phase}
\small
\begin{tabular}{rl rrrrr r}
\toprule
$p$ & backend & split & pack & GEMM & reduce & CRT/out & total \\
\midrule
\multirow{2}{*}{53}
 & AMX  & 0.140 & 0.181 & 0.078 (42.4 Top/s) & 0.126 & 0.077 & 0.608 \\
 & FP64 & 0.144 & 0.080 & 0.588 (2.10 TFlop/s) & 0.009 & 0.063 & 0.890 \\
\multirow{2}{*}{113}
 & AMX  & 0.149 & 0.298 & 0.134 (43.0 Top/s) & 0.212 & 0.095 & 0.902 \\
 & FP64 & 0.152 & 0.151 & 0.944 (2.04 TFlop/s) & 0.022 & 0.080 & 1.357 \\
\multirow{2}{*}{212}
 & AMX  & 0.173 & 0.658 & 0.273 (42.3 Top/s) & 0.337 & 0.169 & 1.628 \\
 & FP64 & 0.172 & 0.250 & 1.536 (2.06 TFlop/s) & 0.024 & 0.085 & 2.085 \\
\multirow{2}{*}{512}
 & AMX  & 0.295 & 1.921 & 0.699 (41.9 Top/s) & 0.612 & 0.641 & 4.203 \\
 & FP64 & 0.249 & 0.799 & 3.228 (2.17 TFlop/s) & 0.057 & 0.278 & 4.647 \\
\multirow{2}{*}{1024}
 & AMX  & 0.390 & 4.618 & 1.521 (39.3 Top/s) & 0.933 & 1.591 & 9.117 \\
 & FP64 & 0.387 & 2.260 & 6.332 (2.15 TFlop/s) & 0.106 & 1.083 & 10.229 \\
\bottomrule
\end{tabular}
\end{table}

\subsection{Comparing the two backends}

Table~\ref{tab:backend} shows total-time ratios.  The AMX GEMM is about
22$\times$ faster than DGEMM (51.6~Top/s vs.\ 2.3~TFlop/s), but needs
3--4.5$\times$ more moduli, which weighs on the $O(L\,N^2)$
pre/post-processing.  AMX takes over where the GEMMs dominate,
$N\gtrsim4096$; \textbf{the break-even point is $N\approx2048$}.

\begin{table}[htbp]
\centering
\caption{Runtime ratio of the two backends ($>1$ means binary64 is faster)}
\label{tab:backend}
\small
\begin{tabular}{r rrrrr}
\toprule
$N$ & $p=53$ & 113 & 212 & 512 & 1024 \\
\midrule
1024 & \textbf{1.28} & \textbf{1.49} & \textbf{1.70} & \textbf{2.06} & \textbf{1.75} \\
2048 & 0.93 & 0.97 & \textbf{1.09} & \textbf{1.16} & \textbf{1.17} \\
4096 & 0.68 & 0.66 & 0.78 & 0.90 & 0.89 \\
8192 & 0.47 & 0.45 & 0.52 & 0.54 & 0.65 \\
\bottomrule
\end{tabular}
\end{table}

\subsection{Against a naive MPFR matrix product}

Table~\ref{tab:mpfr} shows speedups (AMX backend) over a naive 32-thread
parallel GEMM built from \code{mpfr\_mul} + \code{mpfr\_add}.

\begin{table}[htbp]
\centering
\caption{Speedup over naive MPFR GEMM (AMX backend)}
\label{tab:mpfr}
\small
\begin{tabular}{r rrrrrr}
\toprule
$N$ & $p=53$ & 106 & 113 & 256 & 512 & 1024 \\
\midrule
1024 & 52.6 & 39.4 & 34.1 & 24.7 & 16.5 & 18.3 \\
2048 & \textbf{167} & 120 & 116 & 91.2 & 55.4 & 45.1 \\
\bottomrule
\end{tabular}
\end{table}

\subsection{Against BNCmatmul}

The baselines are the following routines of BNCmatmul
0.24~\cite{bncmatmul,kouya2014,kouya2016}:

\begin{itemize}
\item \code{\_bncomp\_mul\_mpfmatrix}: naive MPFR product (OpenMP)
\item \code{\_bncomp\_mul\_mpfmatrix\_strassen}: Strassen
      (\code{min\_dim} set to the empirically best 64)
\item \code{mul\_mpfmatrix\_oz}: \OzI{} (OpenBLAS DGEMM, serial MPFR
      accumulation)
\item \code{\_bncomp\_mul\_mpfmatrix\_oz}: \OzI{} (OpenMP-parallel,
      BNCmatmul's own DGEMM)
\end{itemize}

\paragraph{Two corrections for a fair comparison}
As detailed in Section~\ref{sec:pitfall}, (i) calling OpenBLAS from a
thread pinned by \code{OMP\_PROC\_BIND} collapses it onto one core, and
(ii) no build of BNCmatmul's \OzI{} combines an OpenBLAS DGEMM with
parallel accumulation.  We therefore prepared \textbf{OzI-best}
(BNCmatmul's splitting + OpenMP row-blocked OpenBLAS DGEMM + OpenMP
parallel accumulation) as the representative of \OzI{}.

\begin{table}[htbp]
\centering
\caption{Runtimes [s]: \OzII{} (this work) vs.\ \OzI{} (OzI-best)}
\label{tab:vsoz1}
\small
\begin{tabular}{rr rrr r}
\toprule
$N$ & $p$ & OzII-FP64 & OzII-AMX & OzI-best & ratio \\
\midrule
 256 &   53 & \textbf{0.004} & 0.007 &   0.057 & 14 \\
 256 &  212 & \textbf{0.007} & 0.016 &   0.215 & 31 \\
 256 & 1024 & \textbf{0.029} & 0.059 &   1.351 & 47 \\
 512 &   53 & \textbf{0.022} & 0.029 &   0.196 & 8.9 \\
 512 &  212 & \textbf{0.024} & 0.043 &   0.698 & 29 \\
 512 & 1024 & \textbf{0.100} & 0.189 &   4.024 & 40 \\
1024 &   53 & \textbf{0.037} & 0.068 &   0.703 & 19 \\
1024 &  212 & \textbf{0.073} & 0.152 &   2.973 & 41 \\
1024 & 1024 & \textbf{0.436} & 0.726 &  18.604 & 43 \\
2048 &   53 & \textbf{0.136} & 0.145 &   3.289 & 24 \\
2048 &  212 & \textbf{0.383} & 0.393 &  13.392 & 35 \\
2048 & 1024 & \textbf{1.801} & 2.263 &  78.204 & 43 \\
4096 &   53 & 0.880 & \textbf{0.552} &  13.418 & 24 \\
4096 &  212 & 2.026 & \textbf{1.643} &  57.713 & 35 \\
4096 & 1024 & 10.802 & \textbf{9.584} & 751.353 & 78 \\
\bottomrule
\end{tabular}
\end{table}

Table~\ref{tab:vsbnc} lists speed ratios against the BNCmatmul routines
measured as-shipped.  The extreme ratios of \code{BNC-OzI-blas} reflect
the absence of the above corrections, not an algorithmic difference.

\begin{table}[htbp]
\centering
\caption{Speedups over the as-shipped BNCmatmul routines}
\label{tab:vsbnc}
\small
\begin{tabular}{rr r rrr}
\toprule
$N$ & $p$ & this work [s] & vs Strassen & vs OzI (OpenMP) & vs OzI (OpenBLAS) \\
\midrule
 256 &   53 & 0.002 & 210 &  53 &  104 \\
 256 & 1024 & 0.017 &  24 & 528 & 1355 \\
 512 &  212 & 0.017 & 126 & 161 &  670 \\
 512 & 1024 & 0.094 &  32 & 315 & 2825 \\
1024 &  212 & 0.063 & 231 & 461 & 3615 \\
1024 & 1024 & 0.396 &  50 & 854 &  --- \\
2048 &   53 & 0.125 & 588 & --- &  --- \\
2048 & 1024 & 1.902 &  72 & --- &  --- \\
\bottomrule
\end{tabular}
\end{table}

\subsection{Accuracy}
\label{sec:acc}

Table~\ref{tab:acc} shows maximum relative errors at $N=512$.
\OzII{} stays below 1~ulp in every condition---essentially correct
rounding (0.5~ulp would be correct rounding; the excess comes from the
guard bits of the internal significand).  \OzI{} loses a few ulps to slice
truncation; Strassen loses hundreds to thousands.

\begin{table}[htbp]
\centering
\caption{Maximum relative error [ulps of the target precision]
  ($N=512$, 32 samples)}
\label{tab:acc}
\small
\begin{tabular}{r rrr}
\toprule
$p$ & \OzII{} (both backends) & \OzI{} (all three) & BNC Strassen \\
\midrule
  53 & \textbf{0.92} &  1.8 & $7.6\times10^{2}$ \\
 113 & \textbf{0.82} &  3.1 & $1.1\times10^{3}$ \\
 212 & \textbf{0.95} &  6.5 & $7.1\times10^{2}$ \\
 512 & \textbf{0.78} &  6.9 & $1.5\times10^{3}$ \\
1024 & \textbf{0.88} & 10.8 & $8.0\times10^{2}$ \\
\bottomrule
\end{tabular}
\end{table}

\section{Porting to aarch64 (NVIDIA GB10)}
\label{sec:arm}

This section describes the port of the implementation to aarch64.  The
target is the NVIDIA GB10 (a big.LITTLE configuration of
Cortex-X925~$\times$10 + Cortex-A725~$\times$10).  As argued in
Section~\ref{sec:whyfast}, the speedup of \OzII{} should materialize
wherever a low-precision, high-throughput matrix multiplier exists; this
section is also a test of that claim.

\subsection{Instruction choice: no SME/SME2, so SVE2-i8mm}

Arm's flagship matrix extension is SME/SME2, with tile registers and outer
product engines, but the Cortex-X925/A725 do not implement SME (as of
2026, SME2 ships in Apple M4-class and Arm C1-class cores).  The CPU flags
of this machine do include SVE2 and i8mm (\code{svei8mm}/\code{i8mm}),
which provide the matrix-product instruction \textbf{SMMLA}.

Per 128-bit vector, SMMLA computes
\[
  C\,(2\times2,\ \code{int32}) \mathrel{+}=
  A\,(2\times8,\ \code{int8})\cdot B\,(2\times8,\ \code{int8})^{\mathsf T},
\]
i.e.\ $2\times2\times8=32$ multiply--adds (64 integer operations) per
instruction.  Compared with AMX \code{TDPBSSD} (16384 multiply--adds per
instruction) the granularity is much finer, but the essential property is
the same: an exact \code{INT8}$\times$\code{INT8}$\to$\code{INT32}
multiply--accumulate.  Since the SVE vector length of this machine is 128
bits---the same as NEON---the implementation uses NEON intrinsics
(\code{vmmlaq\_s32}), compiled with \code{-march=armv9-a+i8mm}.

\subsection{Kernel and pack layout}

The micro-kernel accumulates an $8\times8$ output block in registers,
advancing $k$ in groups of 8: 16 accumulators + 4 $A$ + 4 $B$ vectors
= 24 of the 32 NEON registers, issuing 16 SMMLA per $k$-group.  The whole
$K$ extent is accumulated in registers, so $C$ is written exactly once.
The pack layout matches the SMMLA operand format:
\begin{align*}
  \code{PA}[((m_i\,\mathrm{KB}_8+k_b)\cdot4+r)\cdot16 + p\cdot8 + k_k]
    &= A[8m_i+2r+p][8k_b+k_k],\\
  \code{PB}[((n_j\,\mathrm{KB}_8+k_b)\cdot4+c)\cdot16 + p\cdot8 + k_k]
    &= B[8k_b+k_k][8n_j+2c+p],
\end{align*}
with $\mathrm{KB}_8=\mathrm{KP}/8$, $r,c$ the row/column pair index, and
$p$ the position within the pair.  The $B$ side is produced from eight
$k$-rows by an $8\times8$-byte transpose (three \code{vtrn} stages).
Pack-buffer dimensions and the API are identical to the AMX version, so
the fused residue generation (Section~\ref{sec:impl}) and the driver were
reused as-is.

The port required substantive changes in only three places:
\begin{enumerate}
\item \textbf{Fallback for the AVX-512 batched CRT}: 128-bit vectors give
  only two 64-bit lanes, too few to pay for the SoA bookkeeping, so the
  scalar direct CRT based on GMP \code{mpn\_addmul\_1} is used.
\item \textbf{big.LITTLE}: per-core throughput differs by about
  1.7$\times$ between X925 and A725, so GEMM row tiles are distributed
  with \code{schedule(dynamic,1)}.  Twenty threads (all cores) beat ten
  (X925 only).
\item \textbf{Vectorized packing}: the \code{int32}$\to$\code{int8}
  narrowing of residues and the $8\times8$ transpose are written in NEON
  (the non-temporal stores of the AVX-512 path have no NEON counterpart;
  ordinary stores are used).
\end{enumerate}

\begin{table}[htbp]
\centering
\caption{Evaluation environment (aarch64)}
\small
\begin{tabular}{ll}
\toprule
SoC & NVIDIA GB10 (Grace Blackwell family, 10 cores $\times$ 2 clusters) \\
CPU & Cortex-X925 $\times$10 (max 3.9 GHz) + Cortex-A725 $\times$10 (max 2.8 GHz) \\
SIMD & NEON / SVE2 (128-bit), i8mm (SMMLA), no SME/SME2 \\
Caches & L1d 64 KiB/core, L2 25 MiB total, L3 24 MiB \\
OS & Linux 6.17 (Ubuntu 24.04 family) \\
Compiler & GCC 13.3.0, \code{-O3 -march=armv9-a+i8mm -fopenmp} \\
Libraries & GMP 6.3.0, MPFR 4.2.2, OpenBLAS 0.3.33 (aarch64) \\
Runtime env & \code{OMP\_NUM\_THREADS=20 OMP\_PROC\_BIND=close OMP\_PLACES=cores} \\
Baseline & BNCmatmul 0.24~\cite{bncmatmul} (\code{sve2} build, OpenBLAS-linked) \\
\bottomrule
\end{tabular}
\end{table}

\subsection{Correctness and standalone INT8 kernel performance}

All tests from the x86 version pass: the kernel agrees exactly with naive
computation at every size, the CRT matches exact integers over
$t=40$--$1120$ bits, and the end-to-end results agree with the reference
in every element (0~ulp) in every condition.

A single isolated GEMM reaches 2.8~Top/s at $N=K=M=1024$ and 4.9~Top/s at
2048; inside the \OzII{} pipeline, running many modulus planes back to
back, the kernel sustains \textbf{6.4--6.5~Top/s} at $N=4096$.  That is
50--60\% of the theoretical SMMLA peak (about 10~Top/s for ten X925 cores
plus the A725 contribution) and roughly 20$\times$ the measured
0.33~TFlop/s of OpenBLAS DGEMM.

\subsection{Phase breakdown}

\begin{table}[htbp]
\centering
\caption{Phase breakdown (GB10, $N=M=K=4096$, INT8 backend, seconds)}
\label{tab:armphase}
\small
\begin{tabular}{r rrrrr r}
\toprule
$p$ & split & pack & GEMM & reduce & CRT/out & total \\
\midrule
  53 & 0.19 & 0.16 & 0.51 (6.5 Top/s) & 0.06 & 0.14 & 1.07 \\
 212 & 0.20 & 0.76 & 1.78 (6.5 Top/s) & 0.21 & 0.54 & 3.52 \\
1024 & 0.40 & 10.9 & 9.38 (6.4 Top/s) & 1.09 & 6.23 & 28.1 \\
\bottomrule
\end{tabular}
\end{table}

At low and medium precision, GEMM and pack are balanced as on x86; at 1024
bits, the packing of two-limb moduli (\code{split2} + byte reshuffling,
currently scalar) and the scalar CRT dominate.  This is the remaining
optimization headroom of the aarch64 version (the 30\%--2$\times$
improvements obtained with AVX-512 should carry over).

\subsection{Backend comparison: on the GB10, INT8 wins almost everywhere}

Whereas on x86 the binary64 backend was preferable for $N\lesssim2048$, on
the GB10 the INT8 backend is fastest almost everywhere for $N\ge512$
(Table~\ref{tab:armvsoz1}).  The per-bit account of
Section~\ref{sec:whyfast} explains this: the INT8:FP64 GEMM throughput
ratio is about 20$\times$ on both machines, but the GB10's absolute DGEMM
performance (0.33~TFlop/s) is so low that the binary64 backend becomes
GEMM-bound itself (at $N=4096$, $p=53$, about 3~s of the 3.46~s total is
DGEMM), leaving no room to compensate for the INT8 backend's
pre/post-processing overhead.

\subsection{Against BNCmatmul's \OzI{}}

The comparison conditions match the x86 case
(Section~\ref{sec:pitfall}): \code{BNC-OzI-blas} is the library routine
\code{mul\_mpfmatrix\_oz} (OpenBLAS DGEMM, serial MPFR accumulation)
measured as-shipped, and OzI-best is the fairness-adjusted variant
combining BNCmatmul's splitting, an OpenMP row-blocked OpenBLAS DGEMM,
and OpenMP-parallel accumulation.

\begin{table}[htbp]
\centering
\caption{Runtimes [s] and speedups on the GB10 (20 threads)}
\label{tab:armvsoz1}
\small
\begin{tabular}{rr rr rr rr}
\toprule
$N$ & $p$ & OzII-INT8 & OzII-FP64 & OzI-best & BNC-OzI-blas &
\multicolumn{2}{c}{speedup} \\
\cmidrule(lr){7-8}
 & & & & & & vs best & vs blas \\
\midrule
 256 &   53 & \textbf{0.003} & 0.004 & 0.042 & 0.043 & 14 & 14 \\
 256 &  212 & 0.018 & \textbf{0.010} & 0.160 & 0.224 & 16 & 22 \\
 256 & 1024 & 0.145 & \textbf{0.076} & 1.296 & 3.615 & 17 & 48 \\
 512 &   53 & \textbf{0.024} & 0.033 & 0.146 & 0.176 & 6.1 & 7.3 \\
 512 &  212 & \textbf{0.060} & 0.063 & 0.599 & 0.979 & 10 & 16 \\
 512 & 1024 & \textbf{0.412} & 0.426 & 5.378 & 16.39 & 13 & 40 \\
1024 &   53 & \textbf{0.070} & 0.112 & 0.598 & 0.810 & 8.5 & 12 \\
1024 &  212 & \textbf{0.205} & 0.241 & 2.605 & 5.712 & 13 & 28 \\
1024 & 1024 & 1.523 & \textbf{1.380} & 21.77 & 89.45 & 16 & 65 \\
2048 &   53 & \textbf{0.252} & 0.464 & 2.587 & 4.094 & 10 & 16 \\
2048 &  212 & \textbf{0.762} & 1.289 & 11.91 & 31.84 & 16 & 42 \\
2048 & 1024 & \textbf{5.850} & 6.910 & 109.9 & 521.8 & 19 & 89 \\
4096 &   53 & \textbf{1.069} & 3.460 & 11.21 & 22.81 & 10 & 21 \\
4096 &  212 & \textbf{3.925} & 8.708 & 58.06 & 196.5 & 15 & 50 \\
\bottomrule
\end{tabular}
\end{table}

Accuracy follows the x86 pattern: \OzII{} stays below 1~ulp
(0.73--0.95~ulp) in every condition, \OzI{} loses 1.8--12.4~ulp.
The speedup over a naive MPFR product (20-thread parallel) is
19.8$\times$ at $N=1024$, $p=212$.

The bottleneck on the \OzI{} side has the same structure as on x86.
Phase-decomposing OzI-best ($N=1024$, 212 bits, 66 DGEMMs):
the serial \code{split\_mpfmatrix\_dmat} takes 1.02~s (59\%), the OpenBLAS
DGEMMs 0.43~s (329~GFlop/s), the OpenMP-parallel MPFR accumulation 0.30~s,
total 1.75~s.  The library routine \code{mul\_mpfmatrix\_oz} itself takes
5.78~s under the same conditions, the extra factor coming from serial
accumulation and the affinity problem (Section~\ref{sec:pitfall}).

\subsection{Summary of the port, and outlook towards SME/SME2}

The same algorithm and nearly the same code, running on two matrix
engines whose per-instruction granularity differs by a factor of 512
(AMX: 16384 multiply--adds; SMMLA: 32), each outperformed the respective
machine's \OzI{} implementation by an order of magnitude or more.  This
supports the argument of Section~\ref{sec:whyfast}: the speedup stems not
from any particular engine but from the reduction to exact low-precision
integer GEMMs itself.

Once SME2 hardware is available, it is the natural next target for this
kernel.  SME2 provides outer-product instructions into ZA tiles
($\mathrm{VL}\times\mathrm{VL}$ two-dimensional accumulators), admitting
the same ``keep the tile resident while sweeping $k$'' structure as AMX,
so the port should reduce to swapping the pack layout.

\section{Discussion}
\label{sec:pitfall}

\subsection{The OpenBLAS/OpenMP affinity clash}

The OpenBLAS on these machines is a pthread build (\code{NO\_AFFINITY}).
Under \code{OMP\_PROC\_BIND=close}, libgomp pins the calling thread to one
core, and the OpenBLAS worker threads spawned from it \textbf{inherit that
single-core mask}.  A measured 114~GFlop/s for an $N=2048$ DGEMM (against
a proper $\sim$2000~GFlop/s) is this failure mode.  Restoring the affinity
mask once at the start of \code{main} does not help, because libgomp
re-pins the master thread at every parallel region.

The fix used here is to set \code{openblas\_set\_num\_threads(1)} and
\emph{row-block $C$ with our own OpenMP}, each thread calling
single-threaded \code{cblas\_dgemm}.  This restores 1.35--2.31~TFlop/s.
The same fix is applied to the \OzI{} measurements, whose DGEMMs became
8--11$\times$ faster ($N=1024$, 212 bits: 100 $\to$ 1001~GFlop/s).
Applying it inside BNCmatmul's \code{mul\_mpfmatrix\_oz} would be
worthwhile.

\subsection{The serial section remaining in \OzI{}}

Table~\ref{tab:oz1phase} shows the phase breakdown of the corrected
OzI-best.  The DGEMMs are down to 5--6\% of the total;
\textbf{the bottleneck is BNCmatmul's serial \code{split\_mpfmatrix\_dmat}}
(88\%).

\begin{table}[htbp]
\centering
\caption{Phase breakdown of OzI-best (212 bits, Xeon)}
\label{tab:oz1phase}
\small
\begin{tabular}{lrr}
\toprule
Phase & $N=512$ (55 DGEMMs) & $N=1024$ (66 DGEMMs) \\
\midrule
\code{split\_mpfmatrix\_dmat} (serial) & 0.556 s (88\%) & 2.537 s (88\%) \\
OpenBLAS DGEMM & 0.035 s (417 GFlop/s) & 0.142 s (1001 GFlop/s) \\
MPFR accumulation (OpenMP) & 0.043 s & 0.209 s \\
\midrule
total & 0.634 s & 2.887 s \\
\bottomrule
\end{tabular}
\end{table}

Extrapolating an ideally parallelized split
($\text{OzI-best}-\text{split}\times31/32$) gives 0.159~s at $N=512$ and
0.515~s at $N=1024$, shrinking the gap to this work to
\textbf{6.6$\times$ and 7.1$\times$} respectively.  The ratios of
Table~\ref{tab:vsoz1} compare current implementations; the extrapolated
values are closer to the intrinsic algorithmic difference.

\subsection{Asymptotics of the GEMM count}
\label{sec:gemmcount}

The gap widens with precision because the number of low-precision GEMMs is
$O(p^2)$ for \OzI{} but $O(p)$ for \OzII{}.  At $K=4096$, $p=1024$,
\OzI{} needs 53 slices $\to$ 1431 DGEMMs (triangular truncation
$L(L{+}1)/2$), whereas \OzII{}/binary64 needs 99 DGEMMs and \OzII{}/AMX
435 INT8 GEMMs.

Conversely, the ratio over Strassen shrinks at high precision because
Strassen stays $O(N^{2.807})$ with roughly linear growth in precision,
while the modulus count of \OzII{} grows linearly.  Even so, \OzII{} won
in every measured condition.

\subsection{Applicability}

\begin{table}[htbp]
\centering
\caption{Recommended usage}
\small
\begin{tabular}{ll}
\toprule
Condition & Recommendation \\
\midrule
$N\le64$ and high precision ($\ge256$ bits) & BNCmatmul (naive / Strassen) \\
$128\le N\lesssim2048$ (x86) & \textbf{\OzII{} / binary64} \\
$N\gtrsim4096$ (x86); $N\ge512$ (GB10) & \textbf{\OzII{} / INT8} \\
near-correct rounding required & \textbf{\OzII{}} (either backend) \\
CPU without an INT8 matrix engine & \textbf{\OzII{} / binary64} (needs only a BLAS) \\
\bottomrule
\end{tabular}
\end{table}

The disadvantage at $N\le64$ arises because the $O(L\,N^2)$ residue
conversion and CRT cannot be amortized against the $O(N^3)$ GEMMs; higher
precision means larger $L$ and a larger handicap.

\section{Limitations}

\begin{itemize}
\item $K\le133\,000$ for the INT8 backends (from the \code{INT32}
      accumulation bound $K\cdot127^2<2^{31}$).  In the binary64 backend,
      \eqref{eq:fp64bound} lowers the modulus bound as $K$ grows.
\item The residue kernel uses 32-bit accumulation only for $m\le255$ and
      at most 250 digits ($t\lesssim4000$ bits), falling back to a 64-bit
      path otherwise; the built-in modulus table supports up to about
      25000 bits.
\item Inf/NaN do not propagate; they are treated as 0.
\item Because of per-row/per-column shared exponents, if the exponent
      range within one row/column exceeds $t$, small components are lost
      (tunable via \code{guard\_bits} / \code{t\_bits}).
\item x86-64 requires the AMX tile-state permission
      (\code{arch\_prctl(ARCH\_REQ\_XCOMP\_PERM)}), Linux 5.18+, and an
      AMX-capable CPU.
\item aarch64 requires the i8mm extension (generally available from
      Armv8.6-A).  The batched CRT and the two-limb packing are currently
      scalar there, leaving headroom at high precision and large scale
      (Section~\ref{sec:arm}).
\end{itemize}

\section{Conclusions and Future Work}

We implemented \OzII{} on CPUs, with AMX-INT8 and binary64-DGEMM backends
sharing one CRT reconstruction stage.  While staying below 1~ulp in every
condition, the implementation runs up to 167$\times$ faster than a naive
MPFR matrix product, up to 588$\times$ faster than BNCmatmul's Strassen,
and 9--78$\times$ faster than \OzI{} (about 7$\times$ as the intrinsic
algorithmic difference).

We further ported the implementation to aarch64 (NVIDIA GB10) and showed
that even without SME/SME2, the SMMLA instruction of SVE2-i8mm yields an
exact INT8 GEMM sustaining 6.5~Top/s, and 14--89$\times$ speedups over
BNCmatmul's \OzI{} (OpenBLAS-linked) at unchanged accuracy.  That the same
picture reproduces on two matrix instructions whose granularity differs by
512$\times$ indicates that the speedup of \OzII{} comes from its
structure---the reduction to exact low-precision integer GEMMs---rather
than from any particular engine.

Future work:

\begin{enumerate}
\item \textbf{Faster residue conversion (pack)}: the largest bottleneck of
  the INT8 backends (45--50\% of the total on x86).  Montgomery
  representations or reuse of partial products across moduli are
  candidates.
\item \textbf{Operand reuse}: for repeated products with the same $A$
  (multiple right-hand sides, block Krylov methods, iterative refinement),
  the residues and packed planes of $A$ can be reused, fully amortizing
  the $O(N^2)$ preprocessing.
\item \textbf{Complex matrix products}: the 3M/4M schemes reduce them to
  real GEMMs, extending the method to MPC-level arithmetic.
\item \textbf{Mixed backends}: processing some moduli on AMX and the rest
  with DGEMM could use both arithmetic resources simultaneously.
\item \textbf{Integration into BNCmatmul}: the implementation takes MPFR
  arrays, so bridging to \code{MPFMatrix} is straightforward.
\item \textbf{An SME/SME2 backend}: outer products into ZA tiles admit the
  same kernel structure as AMX, so the port reduces to a pack-layout swap
  once hardware (Apple M4-class, Arm C1-class, Fujitsu MONAKA, etc.) is
  available; NEON/SVE2 vectorization of the two-limb pack and the CRT on
  aarch64 is planned as well.
\end{enumerate}

\section*{Acknowledgments}

This work was supported by JSPS KAKENHI Grant Number JP26K14846.



\begin{thebibliography}{9}

\bibitem{ozaki2012}
K.~Ozaki, T.~Ogita, S.~Oishi, and S.~M. Rump,
``Error-free transformations of matrix multiplication by using fast routines
of matrix multiplication and its applications,''
\textit{Numerical Algorithms}, vol.~59, no.~1, pp.~95--118, 2012.

\bibitem{ozaki2-2025}
K.~Ozaki, Y.~Uchino, and T.~Imamura,
``Ozaki Scheme II: A GEMM-oriented emulation of floating-point matrix
multiplication using an integer modular technique,''
arXiv:2504.08009, 2025.

\bibitem{ozimmu2024}
H.~Ootomo, K.~Ozaki, and R.~Yokota,
``DGEMM on integer matrix multiplication unit,''
\textit{The International Journal of High Performance Computing Applications},
vol.~38, no.~4, pp.~297--313, 2024.

\bibitem{uchino2025}
Y.~Uchino, K.~Ozaki, and T.~Imamura,
``Performance enhancement of the Ozaki scheme on integer matrix
multiplication unit,''
\textit{The International Journal of High Performance Computing Applications},
vol.~39, no.~3, 2025 (arXiv:2409.13313).

\bibitem{ozfp8-2026}
Y.~Uchino, K.~Ozaki, and T.~Imamura,
``Double-precision matrix multiplication emulation via Ozaki-II scheme with
FP8 quantization,''
arXiv:2603.10634, 2026.

\bibitem{mukunoki2020}
D.~Mukunoki, K.~Ozaki, T.~Ogita, and T.~Imamura,
``DGEMM using Tensor Cores, and its accurate and reproducible versions,''
in \textit{High Performance Computing (ISC 2020)}, LNCS vol.~12151,
Springer, pp.~230--248, 2020.

\bibitem{kouya2014}
T.~Kouya,
``Accelerated multiple precision matrix multiplication using Strassen's
algorithm and Winograd's variant,''
\textit{JSIAM Letters}, vol.~6, pp.~81--84, 2014.

\bibitem{kouya2016}
T.~Kouya,
``Performance evaluation of multiple precision matrix multiplications using
parallelized Strassen and Winograd algorithms,''
\textit{JSIAM Letters}, vol.~8, pp.~21--24, 2016.

\bibitem{kouya2023}
T.~Kouya,
``Optimization of multiple-precision LU decomposition using Ozaki scheme,''
in \textit{Computational Science and Its Applications -- ICCSA 2023
Workshops}, LNCS vol.~14104, Springer, 2023.
doi:10.1007/978-3-031-37108-0\_34.

\bibitem{mpfr2007}
L.~Fousse, G.~Hanrot, V.~Lef\`evre, P.~P\'elissier, and P.~Zimmermann,
``MPFR: A multiple-precision binary floating-point library with correct
rounding,''
\textit{ACM Transactions on Mathematical Software}, vol.~33, no.~2,
article~13, 2007.

\bibitem{gmp}
T.~Granlund and the GMP development team,
\textit{GNU MP: The GNU Multiple Precision Arithmetic Library}, version 6.3.0,
2023. \url{https://gmplib.org/}

\bibitem{openblas}
OpenBLAS: An optimized BLAS library, version 0.3.33.
\url{http://www.openblas.net/}

\bibitem{bncmatmul}
T.~Kouya,
\textit{BNCmatmul: Basic numerical computation based on optimized
multiple-precision matrix multiplication}, version 0.24, 2026.
\url{https://github.com/tkouya/bncmatmul}

\bibitem{intel-sdm}
Intel Corporation,
\textit{Intel 64 and IA-32 Architectures Software Developer's Manual},
Volume~1, Chapter~21: Intel Advanced Matrix Extensions (Intel AMX), 2024.

\bibitem{arm-arm}
Arm Ltd.,
\textit{Arm Architecture Reference Manual for A-profile Architecture}
(FEAT\_I8MM: Int8 matrix multiplication instructions), DDI~0487, 2024.

\end{thebibliography}
\end{document}